\documentclass[12pt,a4paper]{article}

\usepackage{theorem,enumerate}
\usepackage{amsmath,latexsym,amssymb,amsfonts, epsfig}
\usepackage{eucal}
\usepackage{mathrsfs}
\usepackage{array}
\usepackage{epstopdf} 
\usepackage{MnSymbol}
\newcounter{Scounter}
\theorembodyfont{\normalfont\slshape}
\newtheorem{thm}{Theorem}

\newtheorem{definition}[thm]{Definition}

\newtheorem{observation}[thm]{Observation}

\newtheorem{con}[thm]{Conjecture}
\newtheorem{prob}{Problem}

\def\A{{ \mathcal{A}}}

\bfseries\normalfont

\makeatletter
\def\thanks#1{%
   \footnotemark
   \edef\@tempa{\noexpand\noexpand\noexpand\footnotetext[\the\c@footnote]}%
   \toks@\expandafter{\@thanks}%
   \toks\tw@{{#1}}
   \xdef\@thanks{\the\toks@\@tempa\the\toks\tw@}}
\makeatother

\begin{document}

\title{Special Hist-Snarks}

\author{
Arthur Hoffmann-Ostenhof, Thomas Jatschka}

\date{}
\maketitle

\begin{abstract} 
A Hist in a cubic graph $G$ is a spanning tree $T$ which has only vertices of degree three and one. 
A snark with a Hist is called a Hist-snark, see \cite{HO}. We present several computer generated Hist-snarks which form generalizations of the Petersen graph. Moreover, we state some results on Hist-snarks which have been achieved with computer support.
%It is planed that non-computer related results will appear in an updated version of this paper. At the moment, this paper is an appendix to the paper "Snarks with special spanning trees".
\end{abstract}


\section{Introduction}

All considered graphs are looples and finite. For terminology not defined here, we refer to \cite{BM}. A \textit{cycle} is a $2$-regular connected graph. 
An \textit{end-edge} is an edge incident with a vertex of degree $1$. Edge colorings are considered to be proper edge colorings. A \textit{snark} is a cyclically $4$-edge connected cubic graph of girth $5$ admitting no $3$-edge coloring. 
\textit{Hist} is an abbreviation for homeomorphically irreducible spanning tree, see \cite{ABHT}. 
A Hist in a cubic graph is thus a spanning tree without vertices of degree two. We call a snark $G$ a \textit{Hist-snark} if $G$ has a Hist. For examples of Hist-snarks, see Fig.1-Fig.17 (the Hist is illustrated in bold face). For more informations on Hist-snarks, see \cite{HO}.
The following definition is essential.

\begin{definition}\label{d:}
Let $G$ be a simple cubic graph with a Hist $T$. \\
{\bf(i)} An outer cycle of $G$ is a cycle of $G$ such that all of its vertices are leaves of $T$. \\
{\bf(ii)} Let $\{C_1,C_2,...,C_k\}$ be the set of all outer cycles of $G$ with respect to $T$, then we denote by $oc(G,T)=\{|V(C_1)|,|V(C_2)|,...,|V(C_k)|\}$.
\end{definition}

%If we refer to the outer cycles of a Hist-snark, %then we always mean the outer cycles with
%respect to one given Hist of the graph (a %Hist-snark may have several Hists).
Note that $oc(G,T)$ is a multiset.% several elements of $oc(G,T)$ are possibly the same number. 

Every tree with vertices of degree three and one only, and with at least three edges, is called a \textit{1,3-tree}.
The unique $1,3$-tree which has a $3$-valent vertex with distance $i$ to all leaves, is denoted by 
$T_i$, $i \geq 1$. Hence $T_1$ is isomorphic to $K_{1,3}$ (for an illustration of $T_2$, see Fig.\ref{f:p}). Observe that $T_i$ has a stronger property than the well known parity lemma \cite {Z1} implies for the end-edges of graphs $J$ satisfying $d_J(v) \in \{1,3\} \,\,\forall v \in V(J)$:

\begin{observation} 
Let $T_i$, $i \geq 1$ be $3$-edge colored with colors $1,2,3$. Denote the number of end-edges of $T_i$ having color $j \in \{1,2,3\}$ by $s_j$. Then $s_1=s_2=s_3$.
\end{observation}

We mentioned this observation since it could be helpful for constructing new snarks.\\

We call every snark which has $T_i$ as a Hist, a \textit {$T_i$-snark}. Note that $T_i$-snarks have the smallest possible radius which cubic graphs with order $|V(T_i)|$ can have.
The Petersen graph is the smallest snark and a $T_2$-snark, see Fig.\ref{f:p}. 
The smallest cyclically $5$-edge connected snarks apart from the Petersen graph have already $22$ vertices. They are called Loupekine's snarks and they are both surprisingly $T_3$-snarks, see Figure \ref{f:loupi}. 

\begin{observation} 
There are precisely three $T_3$-snarks.
\end{observation}

Apart form the Loupekine's snarks there is a  third $T_3$-snark denoted by $L_3$ which is defined in the appendix. Note that is has a rare property, namely $L_3$ contains two distinct spanning trees isomporhic to $T_3$. In the subsequent section, we introduce special $T_i$-snarks.
% even rotation $T_3$-snarks (Def.\ref{d:rotation}) which we will see later. 

%Moreover, we define a cubic graph having $T_i$ as a Hist by the vertex sets of the outer cycles, see Fig. \ref{uuuu} for an example.

\section{Rotation snarks}

\begin{figure}[htpb] 
\centering\epsfig{file=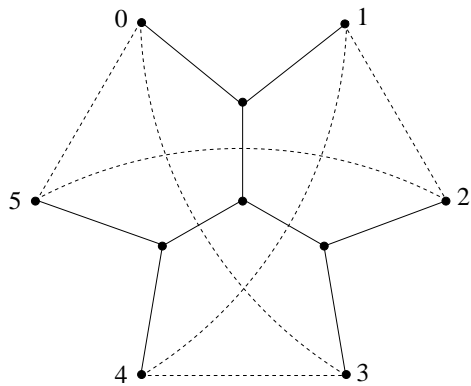,width=2.4in}
\caption{The Petersen graph with a spanning tree isomorphic to $T_2$ illustrated in bold face.}
\label{f:p}
\end{figure}

\begin{figure}[htpb] 
\centering\epsfig{file=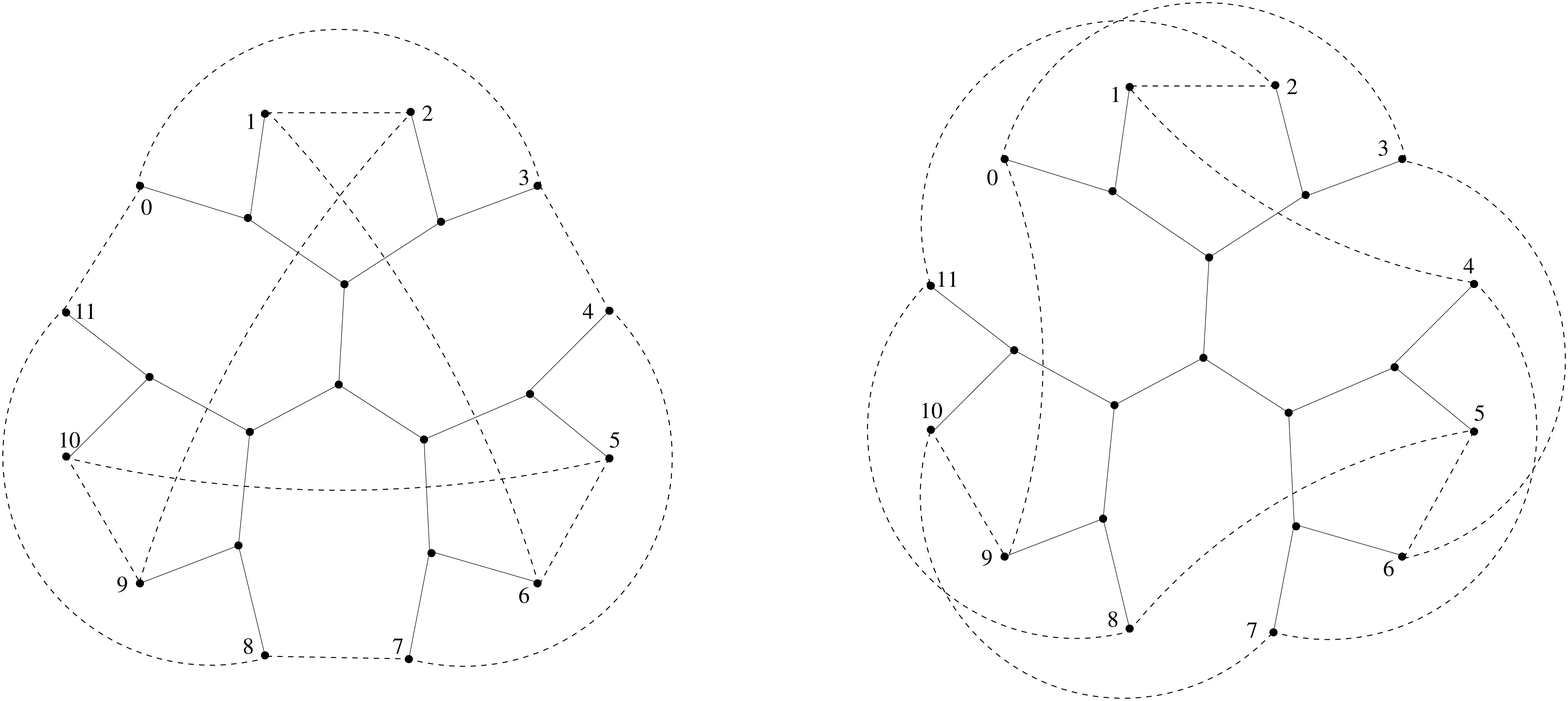,width=6.5in}
\caption{The first Loupekine's snark (on the left side) and the second Loupekine's snark (on the right side).}
\label{f:loupi}
\end{figure}

Considering Fig.\ref{f:p}, we notice that the illustrated $T_2$-snark (the Petersen graph) has a $2\pi/3$ rotation symmetry. We want to generate $T_i$-snarks with this type of symmetry since they form a natural generalization of the Petersen graph and the Loupekine's snarks, see Fig.2. To describe this symmetry, we need some definitions which we state in a short manner. \\A 
\textit{curve} is meant to be a continuous image of a closed unit line segment. A curve is called \textit{simple} if it does not intersect itself.
A \textit{drawing} of a graph $G$ is a representation of $G$ in the plane by representing every vertex $v \in V(G)$ by a distinct point $v'$ in the plane and every edge $xy \in E(G)$ by a simple curve $x'y'$ with endpoints $x'$ and $y'$. Theses edge-curves may cross each other and if they do not cross, then the drawing is called a \textit{planar drawing}. 
We also assume that every vertex-point $v'$ is only part of an edge-curve $x'y'$ if $v' \in \{x',y'\}$.\\ For reasons of simplicity, we call a 
vertex-point a vertex, an edge-curve an edge and we use the same notation for vertex and vertex-point. 

\begin{definition}\label{d:rotation}
Let $G$ be a $T_i$-snark for some $i \geq 2$ and let $l_i$ denote the number of leaves of $T_i$.  
Then $G$ is called a rotation $T_i$-snark (in short a rotation snark) if there is a drawing of $G$ and a labeling of the leaves of the Hist such that all of the following hold:\\
(1) The drawing of $G$ contains a planar drawing of a tree $T$ isomorphic to $T_i$.\\
(1a) All leaf points of $T$ are on an imaginary circle $C$ (we introduce $C$ for defining the leaf labels). All edges of $T$ are within the finite region of $C$. The leaves are labeled from $0$ to $l_i-1$ in cyclic clockwise order with respect to $C$.\\
(1b) Let $H$ be the graph defined by the planar drawing consisting of the already defined drawing of $T$ and the $l_{i}$ arcs of $C$ (regarding them as edges) connecting the leaves of $T$. Then the two 
leaves $l_{i}-1$, $0$ and the central vertex $x \in V(T)$ with $d_T(x,v)=i$ for all leaves $v \in V(T)$, are all together contained in one facial cycle of $H$. \\
(2) Every edge $ab \in E(G)$ which joins two leaves, implies $a+l_i/3 \,\,\,b+l_i/3 \in E(G)$ where $a,b \in \{0,1,2,...,l_i-1\}$ and addition is considered modulo $l_i$. \\
Finally, we denote by $\tilde T_i$ the set of all rotation $T_i$-snarks.
\end{definition}

The Petersen graph is obviously a rotation $T_2$-snark since it holds that $ab \in E(G)$ implies $a+2 \,\,\,b+2 \in E(G)$ for $a,b \in \{0,1,2,...,5\}$, see Fig.\ref{f:p}. For other examples of rotation snarks, see Fig.2-Fig.17.\\
The following definition is very useful and has already been used extensively in \cite{HO}. 

\begin{definition}\label{d:S}
Let $S$ be a multiset of positive numbers, then 
$S^*$ denotes the set of all snarks $G$ which have a Hist $T_G$ such that $oc(G,T_G)=S$.
\end{definition}

%We present the first theorem. 

\begin{thm}\label{t:hh}
(1) The Petersen graph and the two Loupekine's snarks are the only rotation $T_i$-snarks with $i \leq 3$. (2) There are precisely $15$ rotation $T_4$-snarks. In particular,\\ 
$$|\{24\}^* \cap \tilde T_4|=2\,,\,
|\{12,12\}^* \cap \tilde T_4|=1\,,\,
|\{18,6\}^* \cap \tilde T_4|=1\,,\,
|\{8,8,8\}^* \cap \tilde T_4|=8$$ and
$$|\{6,6,6,6\}^* \cap \tilde T_4|=3\,\,.$$
\end{thm}

\bigskip

Below (see Fig.1- Fig.17), are the rotation snarks which are counted in Theorem \ref{t:hh}.

\bigskip

\bigskip

\begin{figure}[htpb] 
\centering\epsfig{file=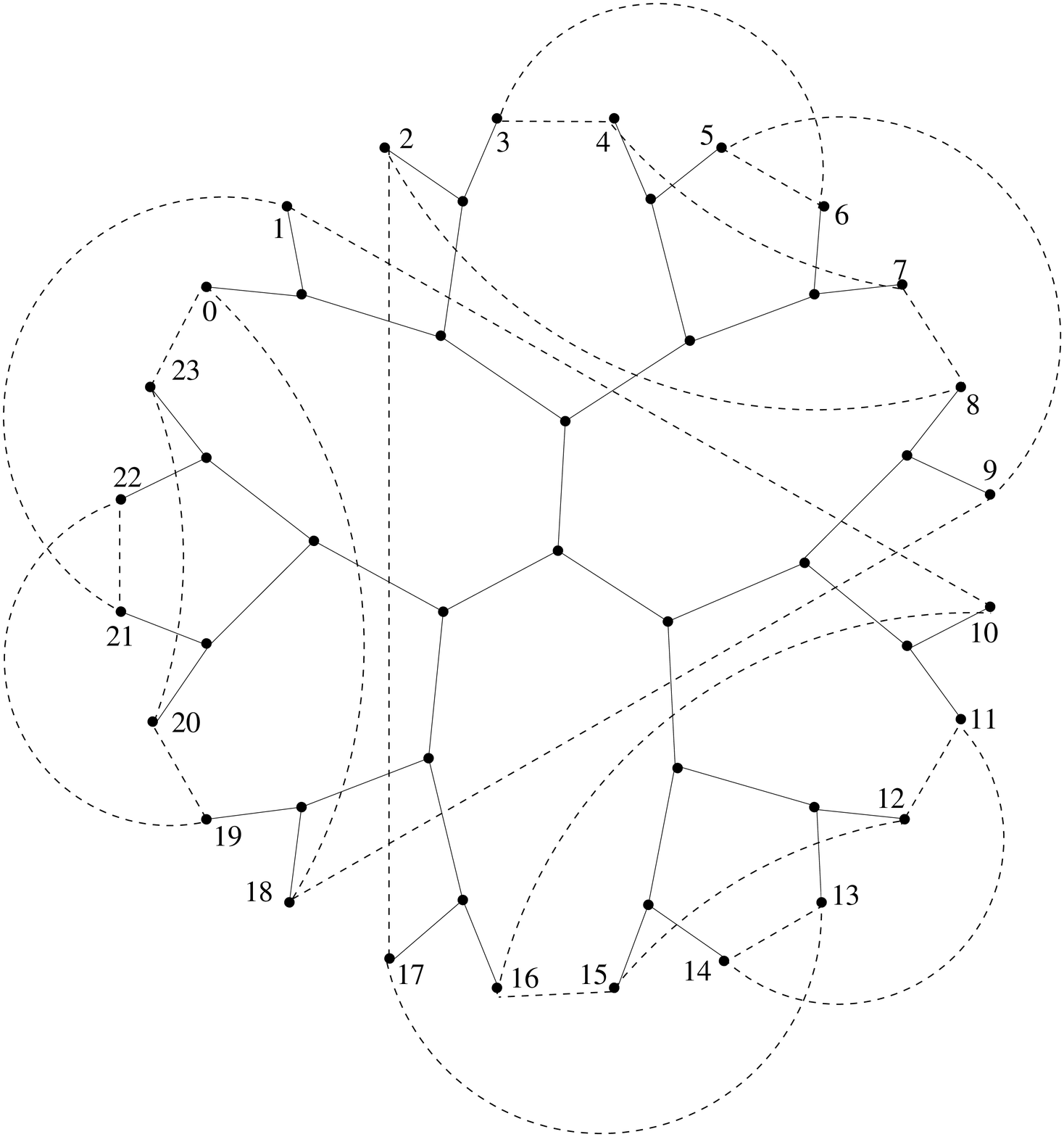,width=4.0in}
\caption{The Hist-snark $H0(24)$.}
\label{f:024}
\end{figure}

\begin{figure}[htpb] 
\centering\epsfig{file=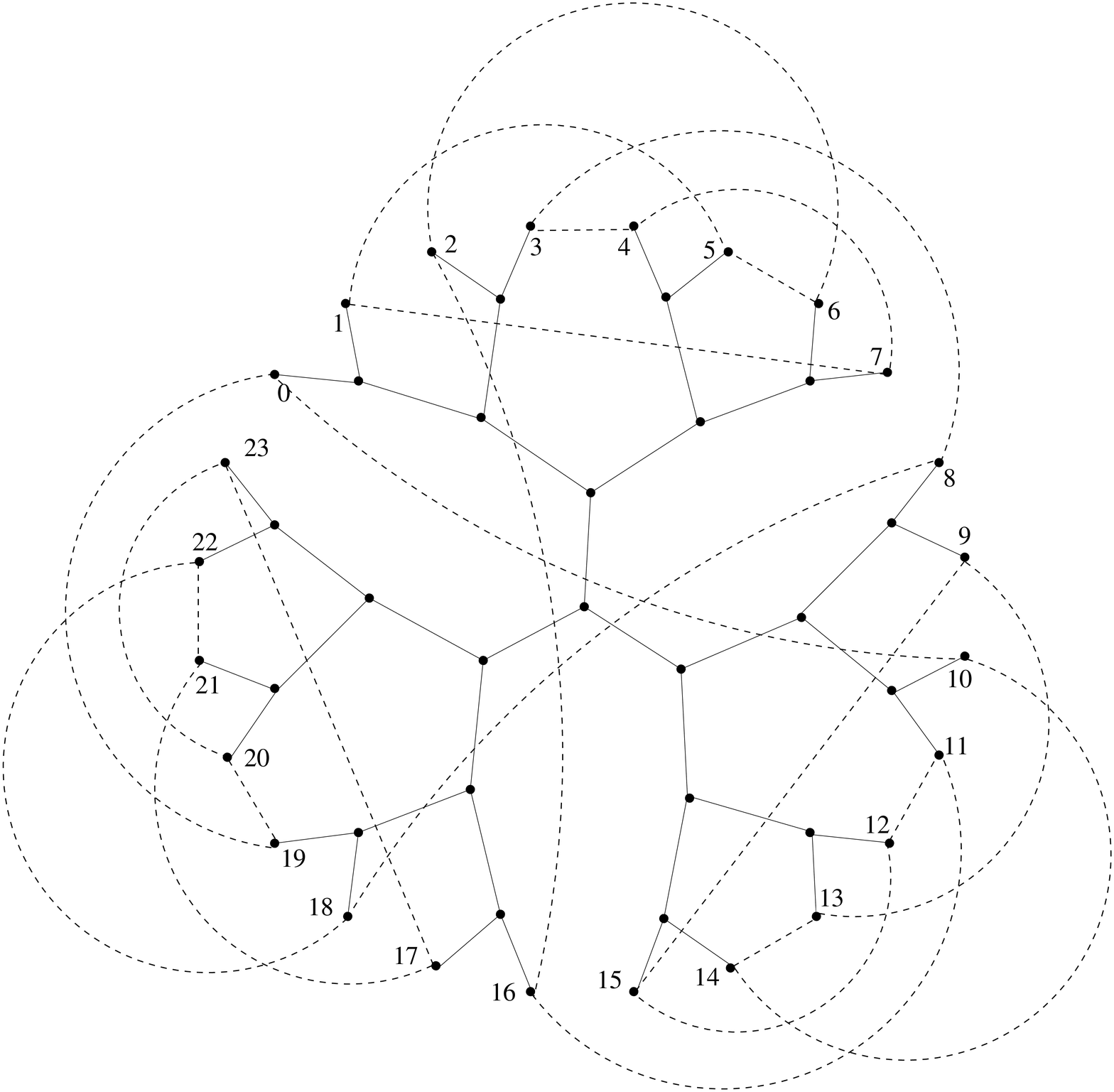,width=4.0in}
\caption{The Hist-snark $H1(24)$.}
\label{f:124}
\end{figure}

\begin{figure}[htpb] 
\centering\epsfig{file=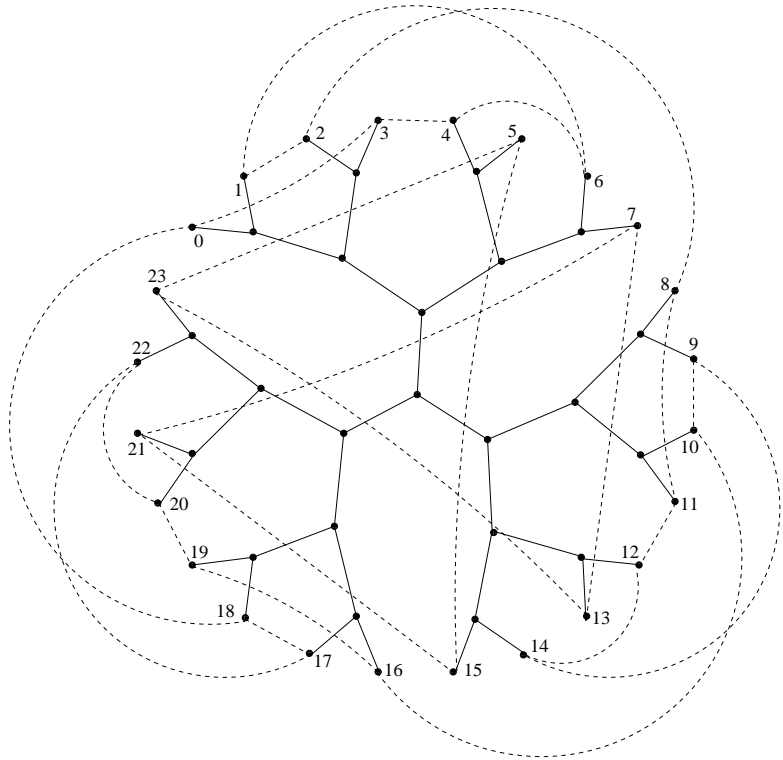,width=4.0in}
\caption{The Hist-snark $H(18,6)$.}
\label{f:186}
\end{figure}

\begin{figure}[htpb] 
\centering\epsfig{file=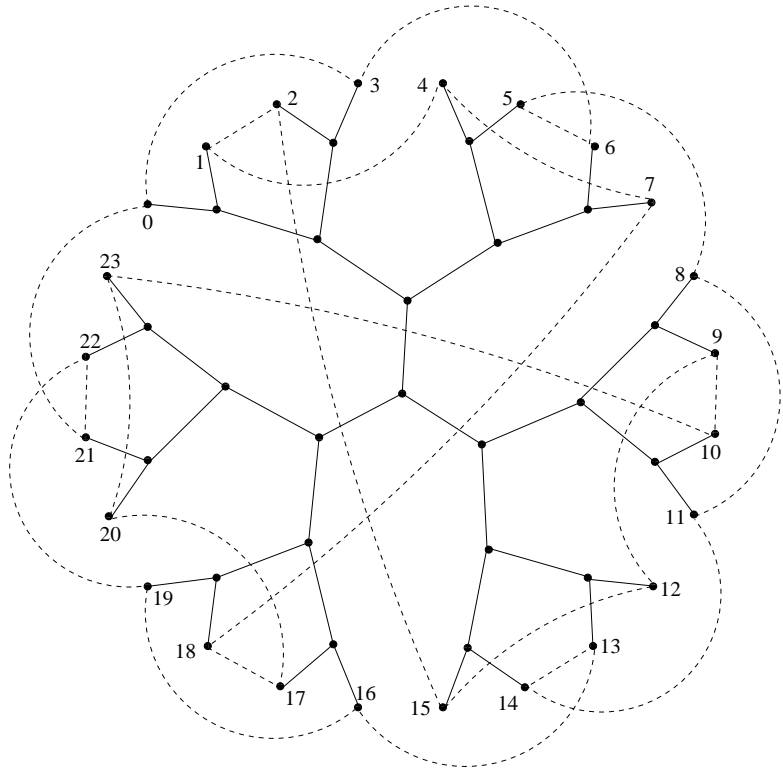,width=4.0in}
\caption{The Hist-snark $H(12,12)$.}
\label{f:1212}
\end{figure}

\begin{figure}[htpb] 
\centering\epsfig{file=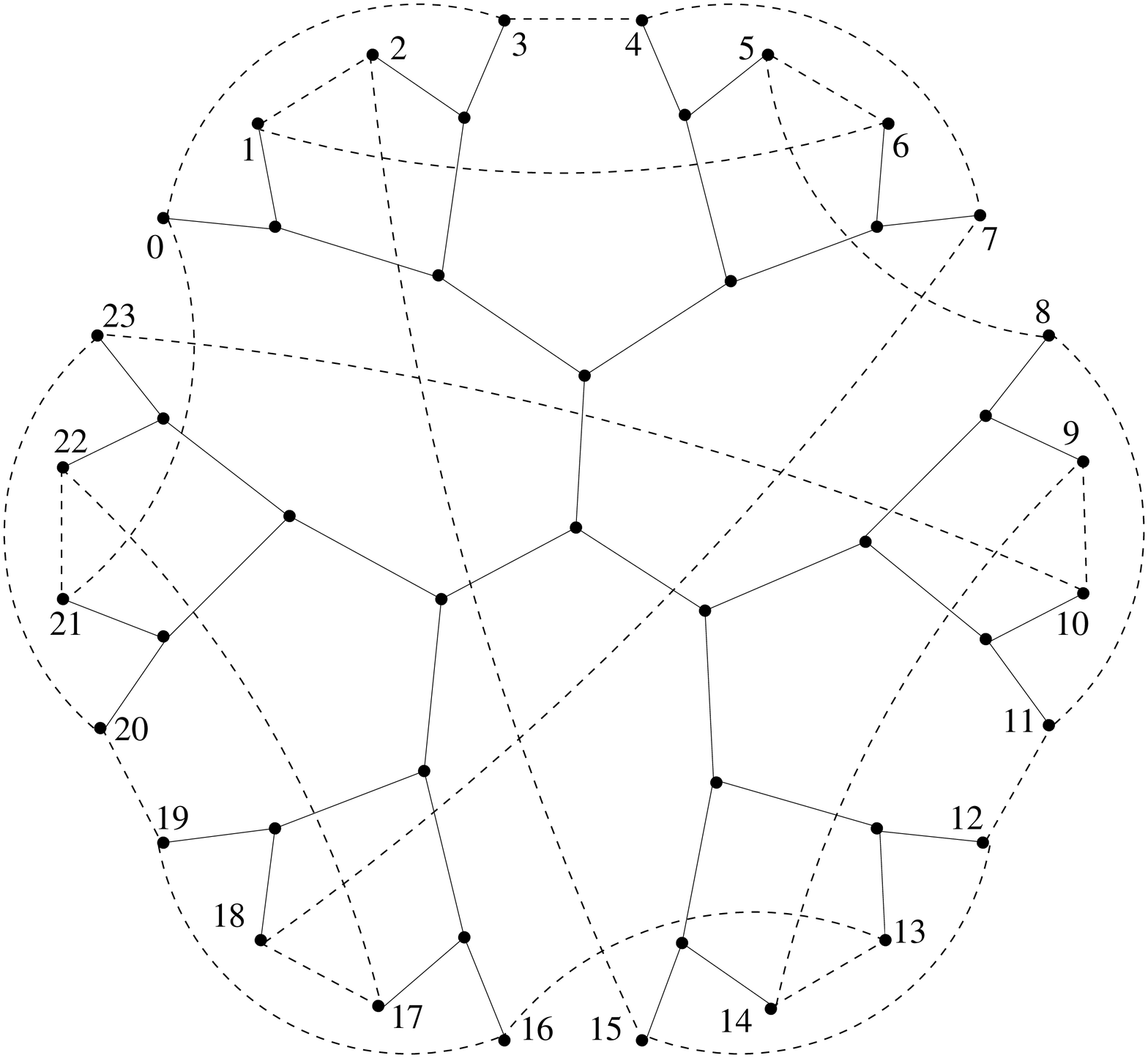,width=4.0in}
\caption{The Hist-snark $H0(8,8,8)$.}
\label{f:0888}
\end{figure}

\begin{figure}[htpb] 
\centering\epsfig{file=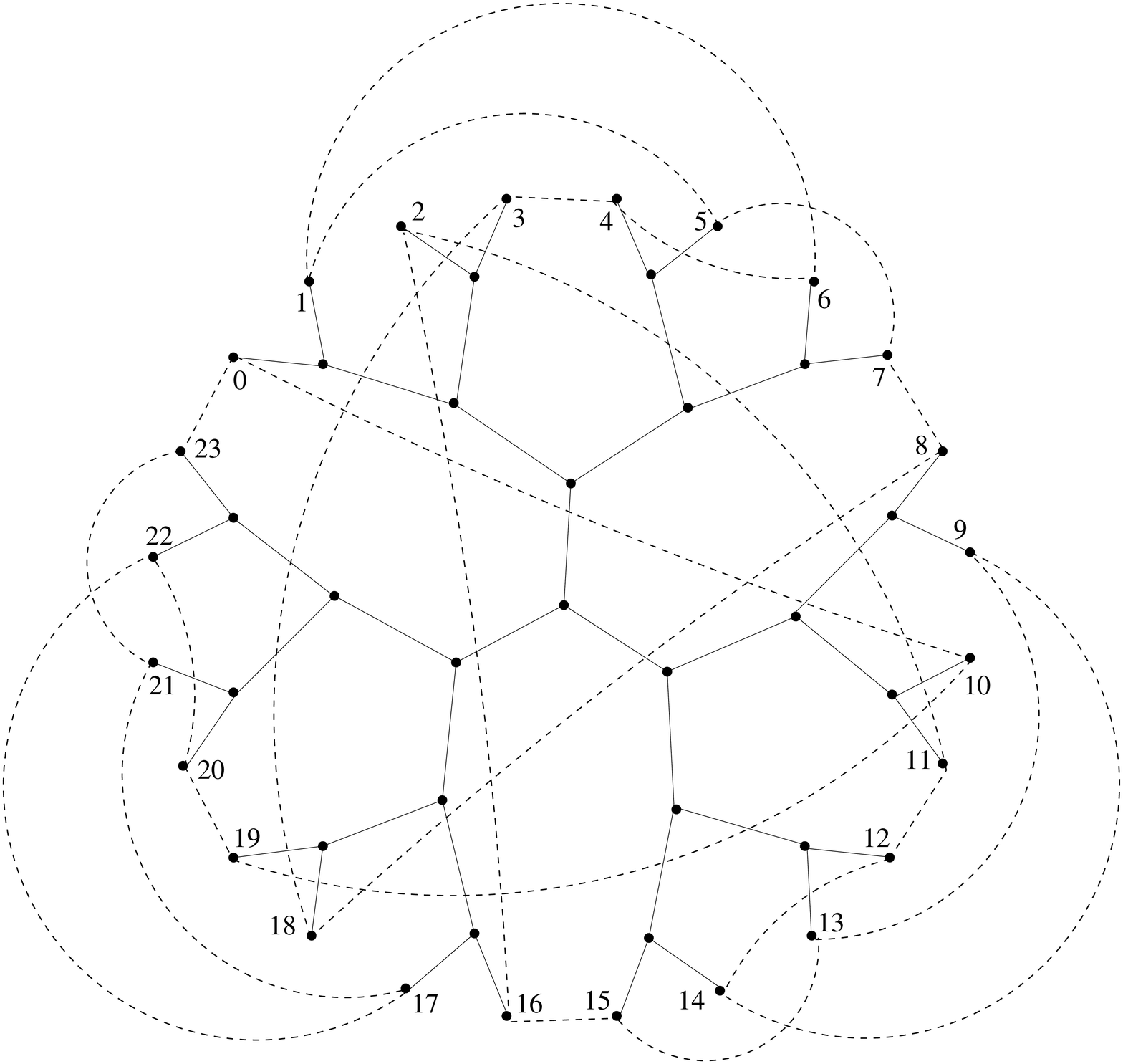,width=4.2in}
\caption{The Hist-snark $H1(8,8,8)$.}
\label{f:1888}
\end{figure}

\begin{figure}[htpb] 
\centering\epsfig{file=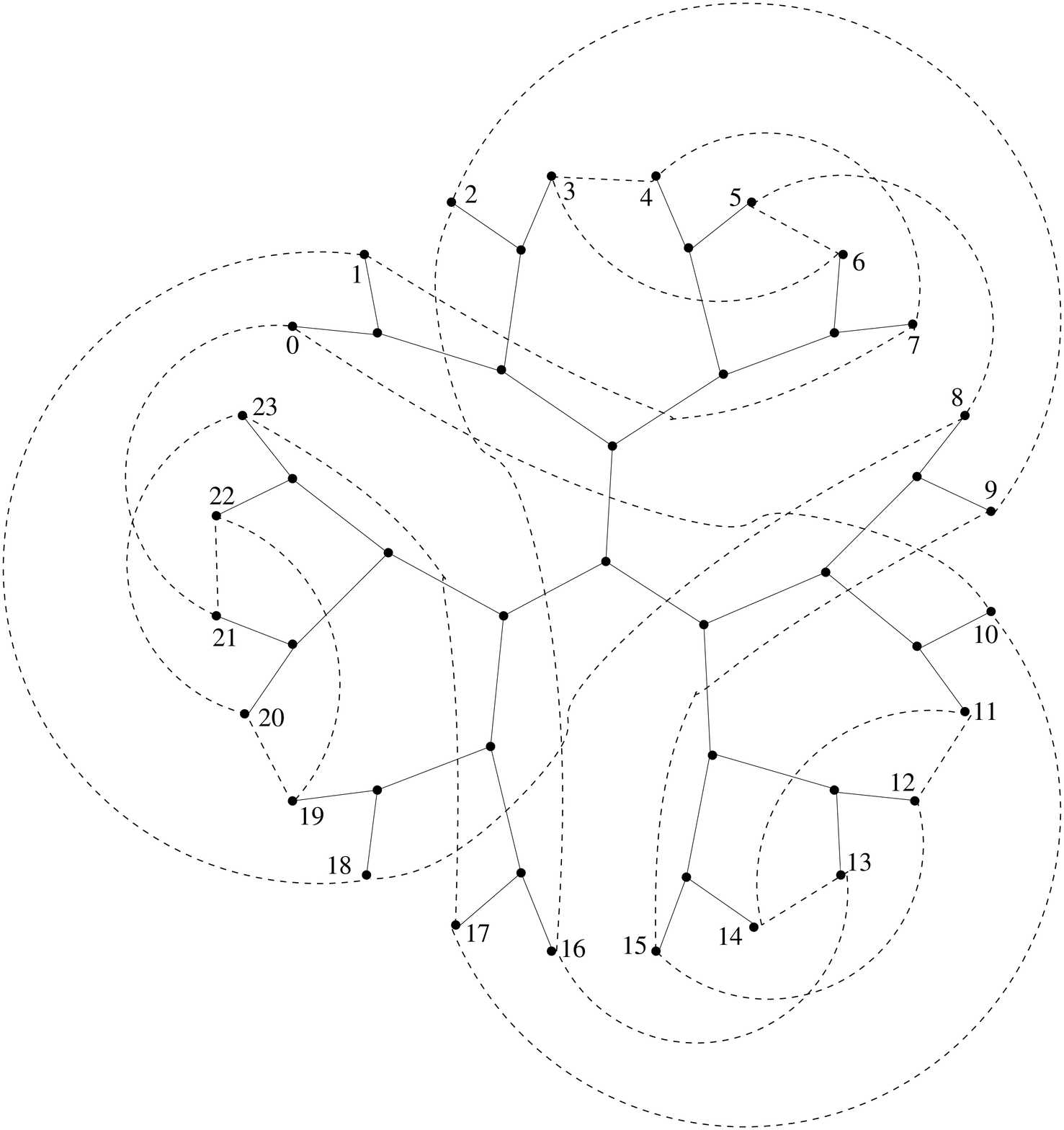,width=4.1in}
\caption{The Hist-snark $H2(8,8,8)$.}
\label{f:2888}
\end{figure}

\begin{figure}[htpb] 
\centering\epsfig{file=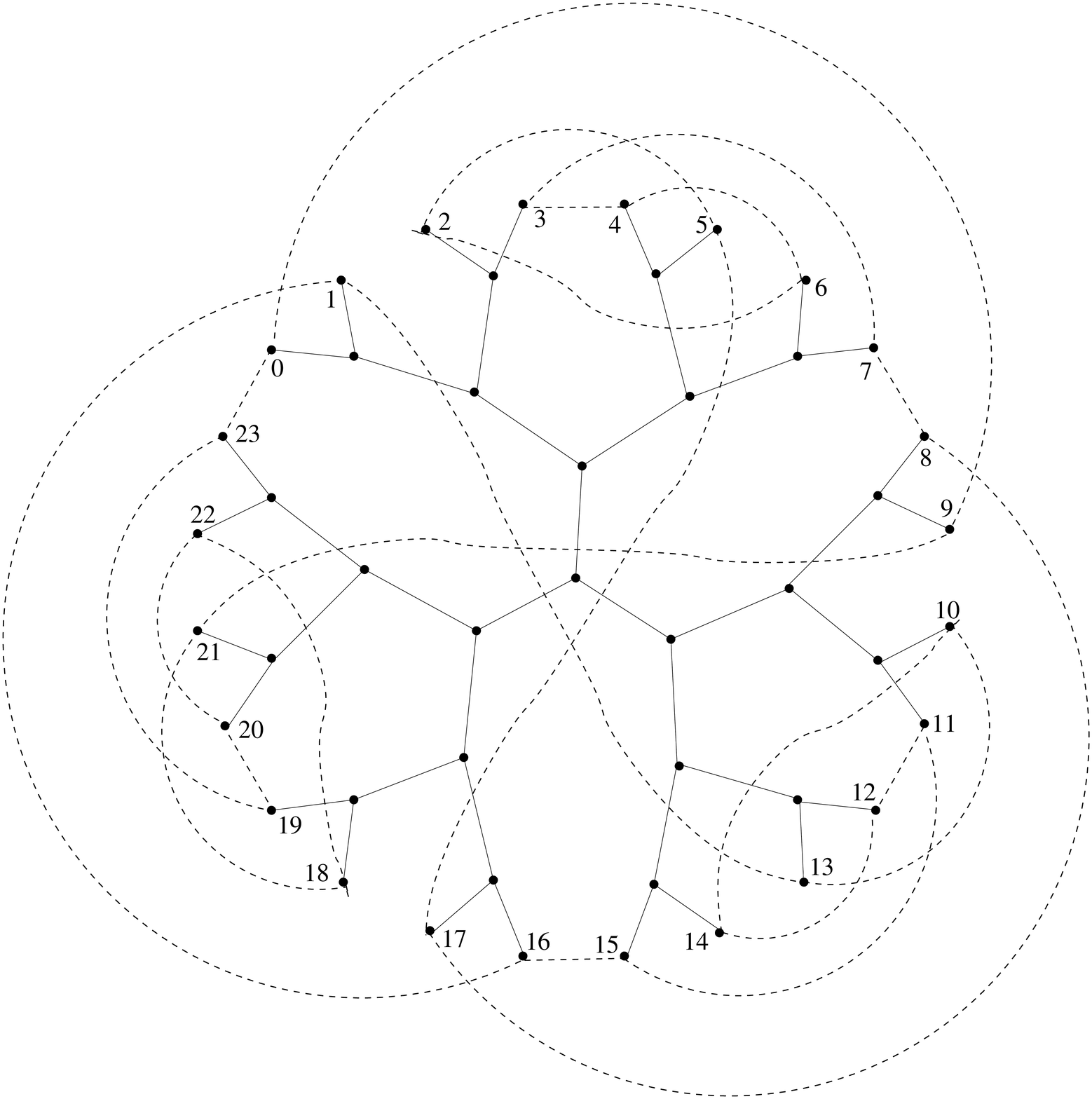,width=4.2in}
\caption{The Hist-snark $H3(8,8,8)$.}
\label{f:3888}
\end{figure}

\begin{figure}[htpb] 
\centering\epsfig{file=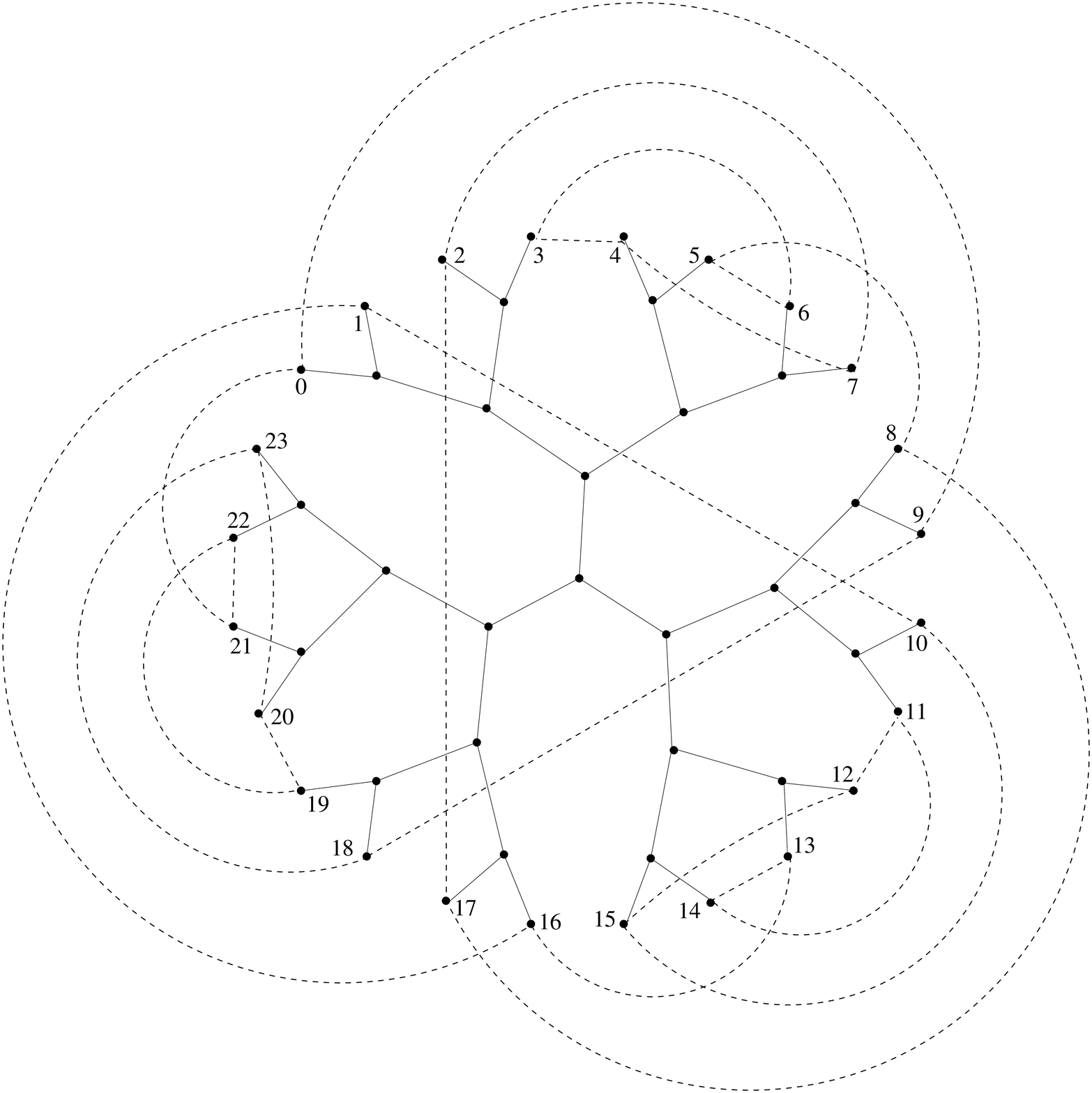,width=4.2in}
\caption{The Hist-snark $H4(8,8,8)$.}
\label{f:4888}
\end{figure}

\begin{figure}[htpb] 
\centering\epsfig{file=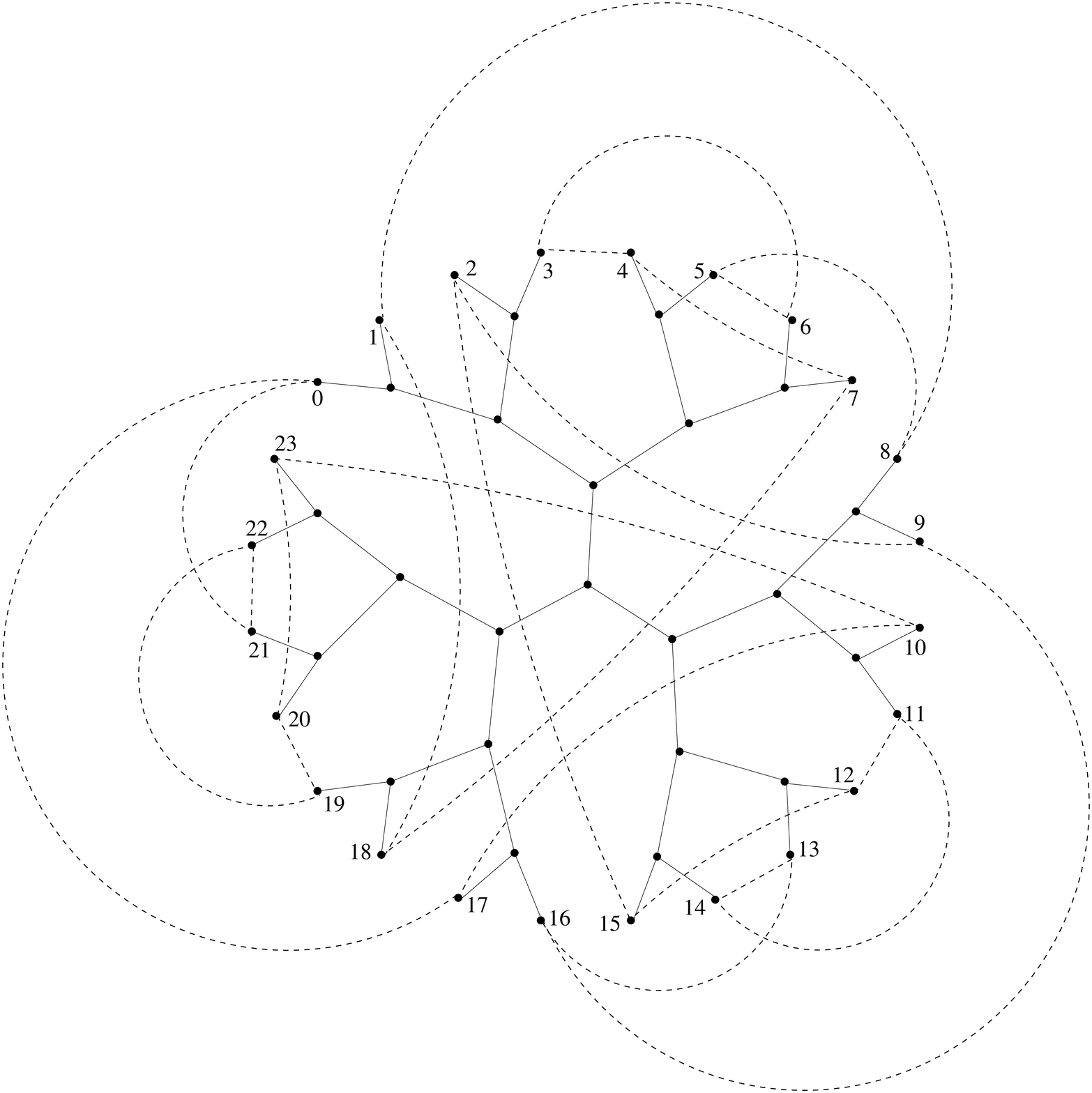,width=4.2in}
\caption{The Hist-snark $H5(8,8,8)$.}
\label{f:5888}
\end{figure}

\begin{figure}[htpb] 
\centering\epsfig{file=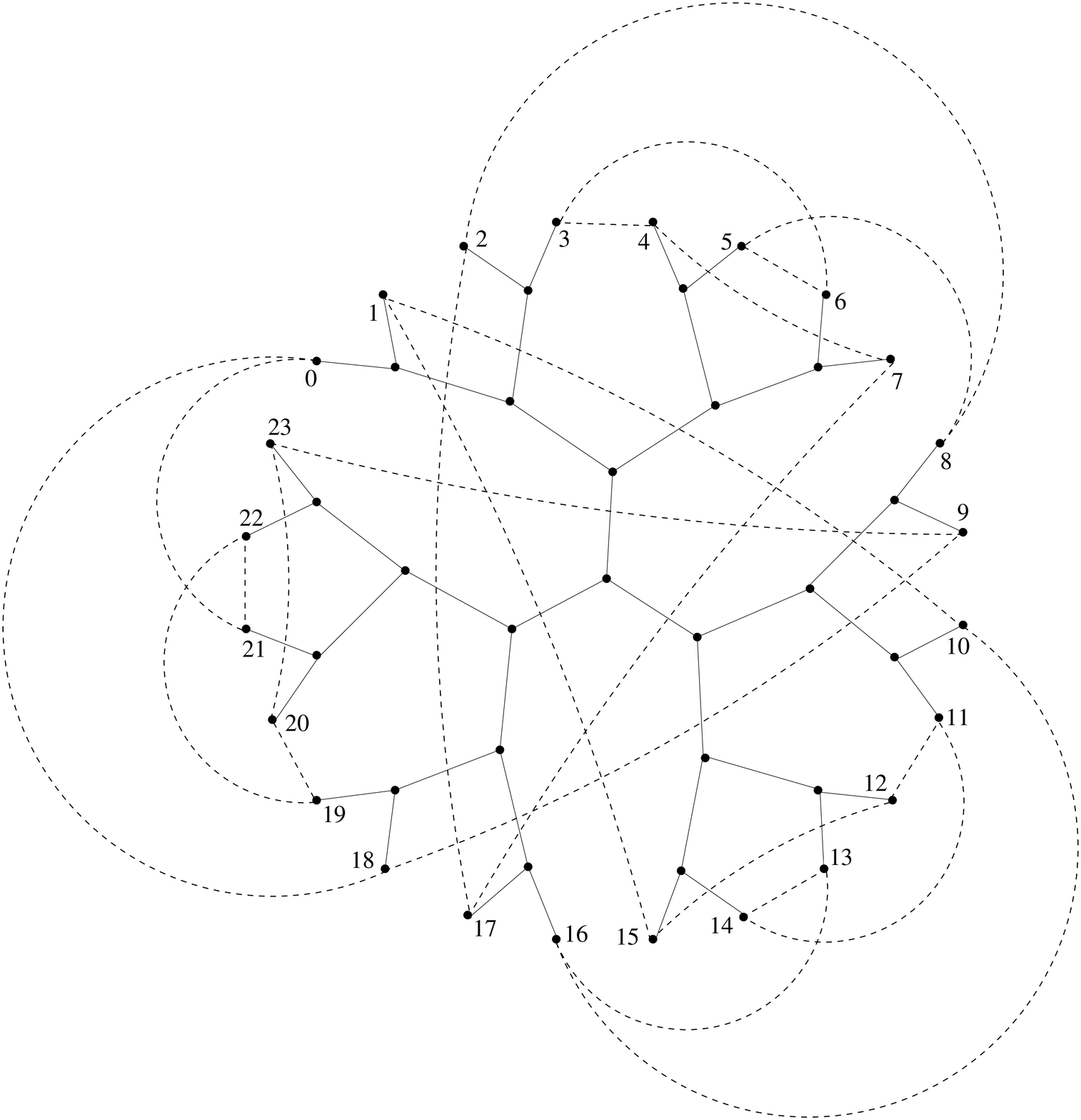,width=4.2in}
\caption{The Hist-snark $H6(8,8,8)$.}
\label{f:6888}
\end{figure}

\begin{figure}[htpb] 
\centering\epsfig{file=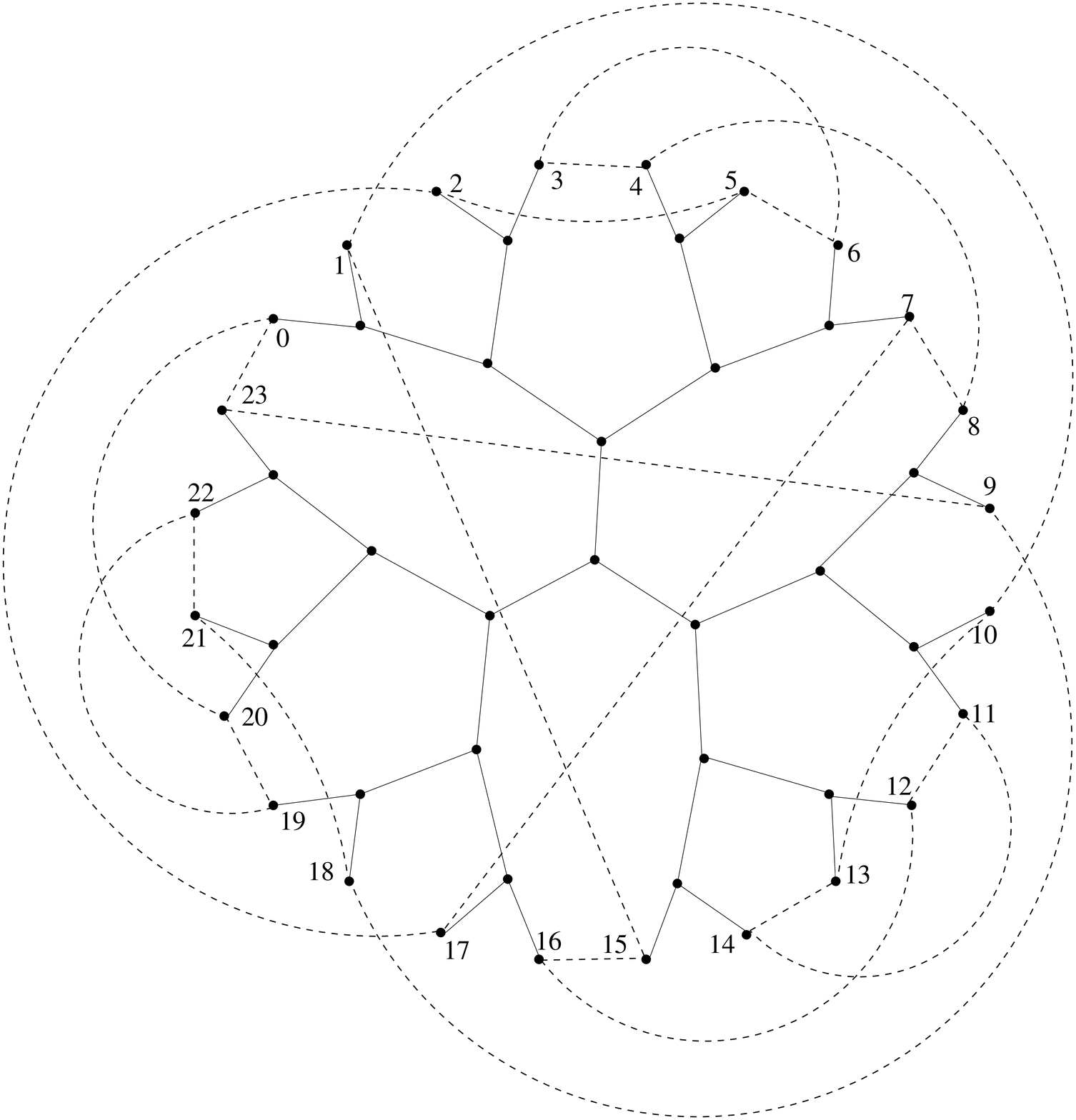,width=4.2in}
\caption{The Hist-snark $H7(8,8,8)$.}
\label{f:7888}
\end{figure}

\begin{figure}[htpb] 
\centering\epsfig{file=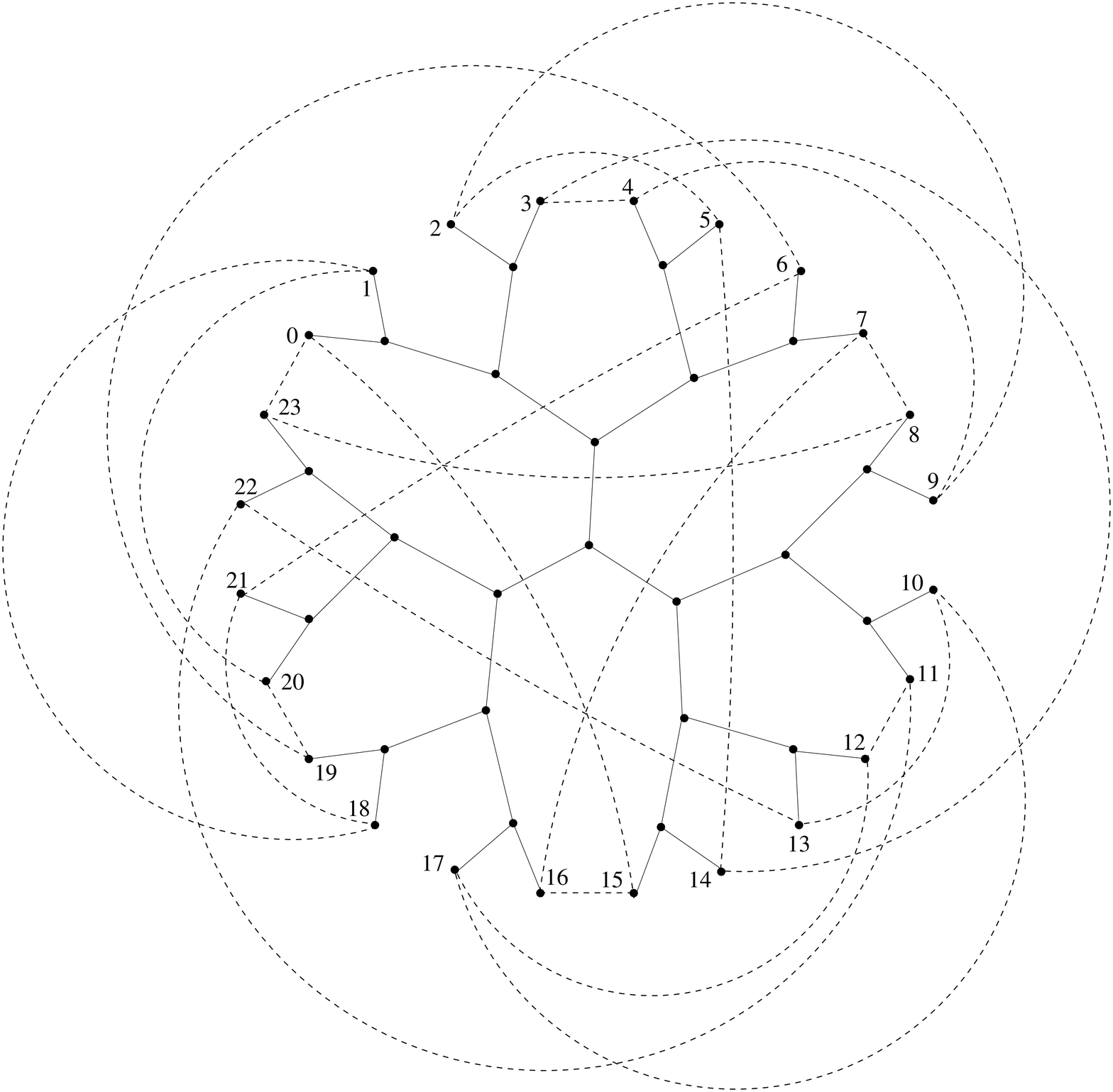,width=4.2in}
\caption{The Hist-snark $H0(6,6,6,6)$.}
\label{f:06666}
\end{figure}

\begin{figure}[htpb] 
\centering\epsfig{file=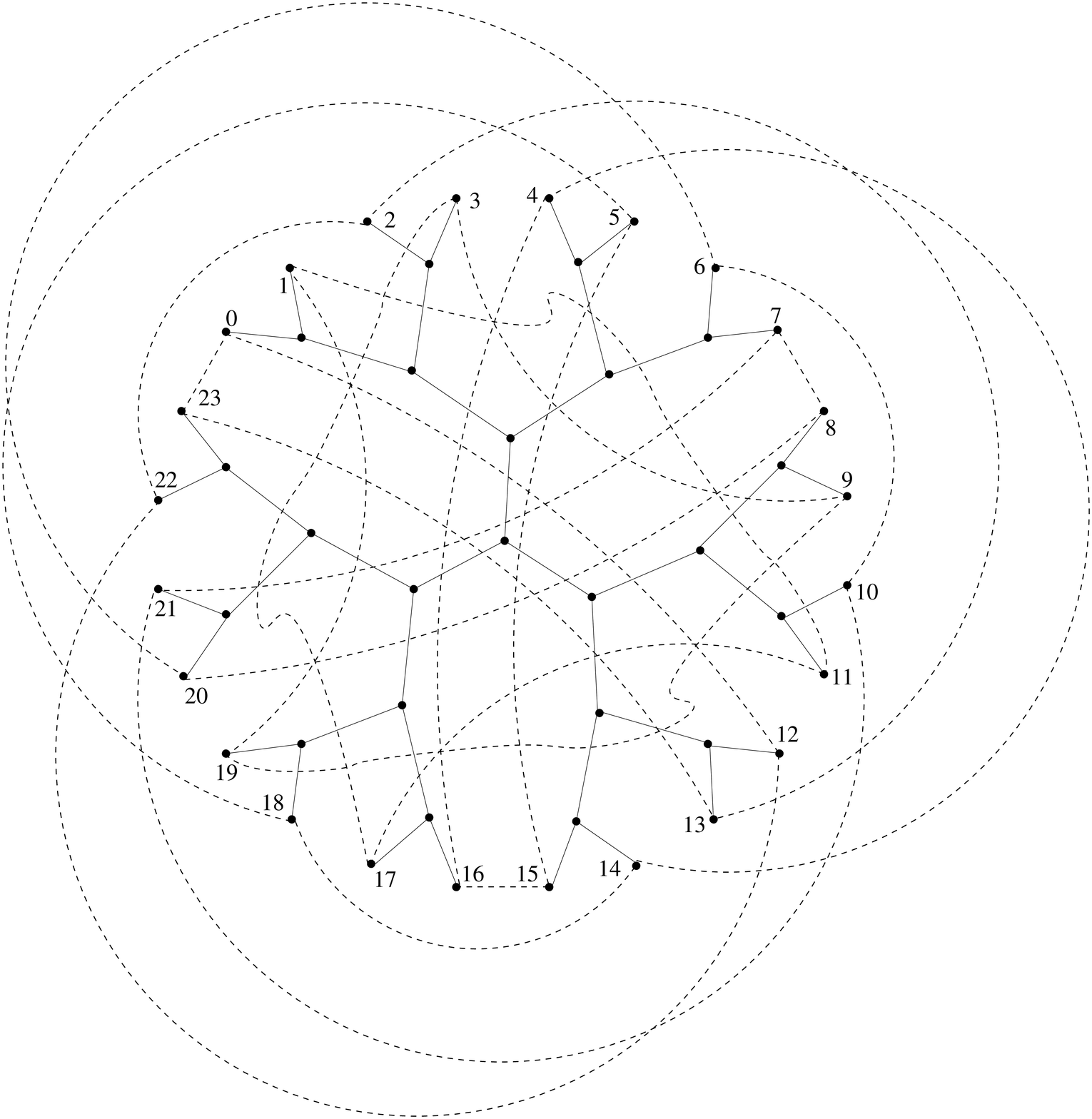,width=4.1in}
\caption{The Hist-snark $H1(6,6,6,6)$.}
\label{f:16666}
\end{figure}

\begin{figure}[htpb] 
\centering\epsfig{file=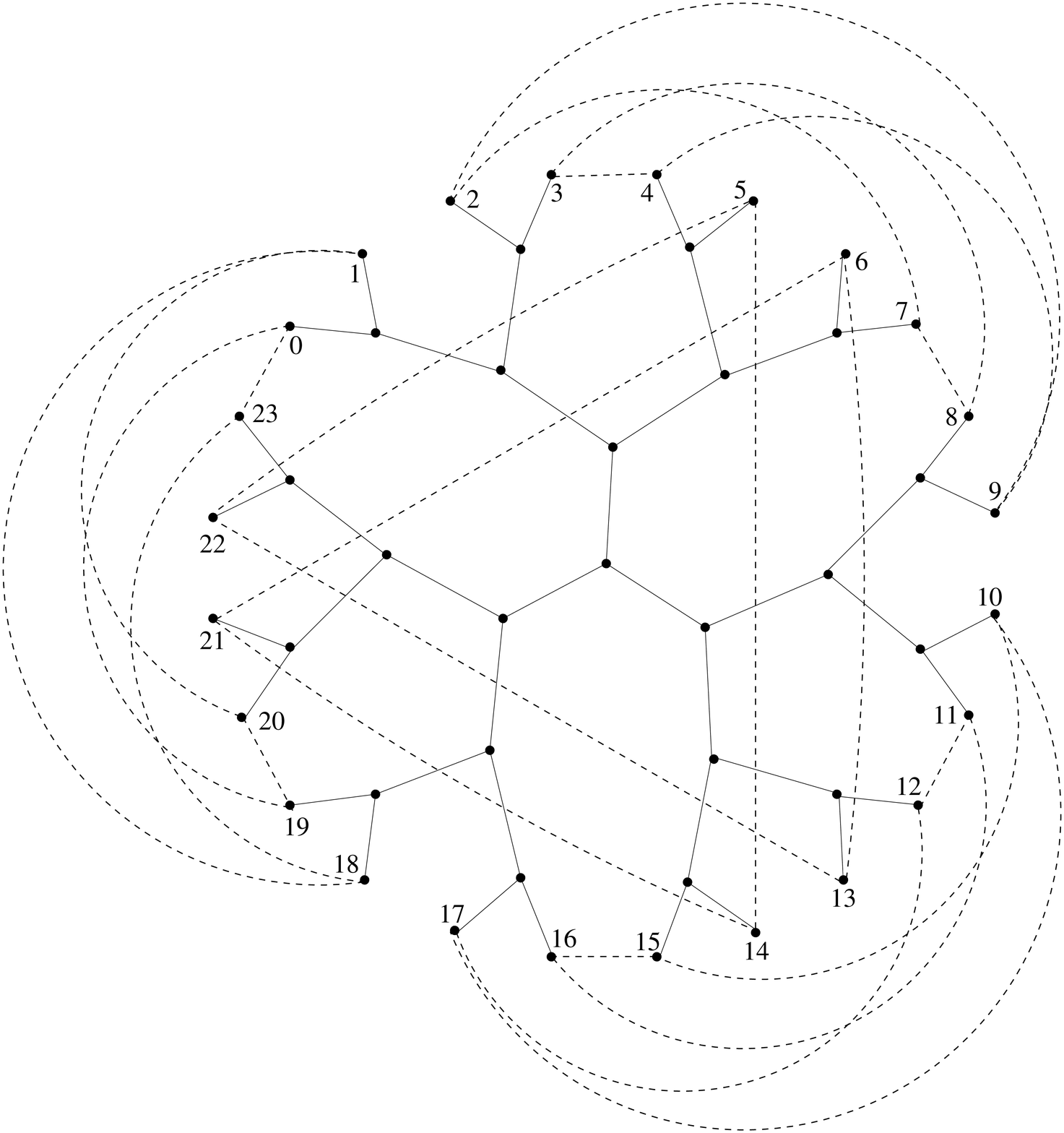,width=4.1in}
\caption{The Hist-snark $H2(6,6,6,6)$.}
\label{f:26666}
\end{figure}

\section{Other results on $T_4$-snarks}

\begin{thm}\label{t:girth}
Let $G$ be cyclically $4$-edge connected cubic graph with girth at least $6$. 
Then $G$ is $3$-edge colorable if $G$ has a spanning tree isomorphic to $T_i$ for some 
$i \in \{1,2,3,4\}$. 
\end{thm}

\begin{thm}\label{t:aut}
Let $G$ be a $T_4$-snark satisfying $G \in \{6,6,6,6\}^*$. Then the number of automorphisms of $G$ is at most 
$128$. Moreover, the $T_4$-snark $Y \in \{6,6,6,6\}^*$ in the appendix below has precisely $128$ automorphisms.
\end{thm}

\section{Open problems}

One motivation to consider $T_i$-snarks was the intention to obtain sooner or later small snarks of girth $7$ or even cyclically $7$-edge connected snarks. Since the generated snarks do not even have girth $6$, we state some less ambitious open problems.

\begin{prob}
Is there a $T_i$-snark with girth $k \geq 6$ for some $i$?
\end{prob}

\begin{prob}
Is there a cyclically $k$-edge connected $T_i$-snark with $k \geq 6$ for some $i$?
\end{prob}

\begin{prob}
Is there a cyclically $k$-edge connected rotation $T_i$-snark with $k \geq 6$ for some $i$?
\end{prob}

\begin{con}
There are infinitely many rotation snarks.
\end{con}

\section*{Acknowledgments}
A.Hoffmann-Ostenhof was supported by the Austrian Science Fund (FWF) project P 26686.
The computational results presented have been achieved using the 
Vienna Scientific Cluster.

\section{Appendix}
For the sake of completeness, all above drawn Hist-snarks are presented now in a short manner, namely via the outer cycles whose vertices are presented within brackets in cyclic order. Moreover, at the end of the list, the $T_4$-snark $Y$ is defined, see Theorem \ref{t:aut}. Note that $Y$ is not a rotation snark.\\

\noindent

%$T(8,8,8):=[0,3,4,7,18,17,22,21]\,\,[1,2,15,12,11,8,5,6]\,\,[9,10,23,20,19,16,13,14]$

$H0(24):=[0,18,9,5,6,3,4,7,8,2,17,13,14,11,12,15,16,10,1,21,22,19,20,23]$

$H1(24):=[0,19,20,23,17,21,22,18,8,3,4,7,1,5,6,2,16,11,12,15,9,13,14,10]$

$H(12,12):=[0,21,22,19,16,13,14,11,8,5,6,3]\,\,[ 10,23,20,17,18,7,4,1,2,15,12,9]$

$H(18,6):=[0,18,17,22,20,19,16,10,9,14,12,11,8,2,1,6,4,3 ]\,\, [13,23,5,15,21,7]$
 
$H0(8,8,8):=[0,3,4,7,18,17,22,21 ]\,\,[1,2,15,12,11,8,5,6 ]\,\,[9,10,23,20,19,16,13,14]$

$H1(8,8,8):=[0,23,21,17,22,20,19,10 ]\,\,[ 8,7,5,1,6,4,3,18 ]\,\,[ 16,15,13,9,14,12,11,2]$

$H2(8,8,8):=[0,21,22,19,20,23,17,10 ]\,\,[ 8,5,6,3,4,7,1,18 ]\,\,[ 16,13,14,11,12,15,9,2]$

$H3(8,8,8):=[0,23,19,20,22,18,21,9 ]\,\,[ 8,7,3,4,6,2,5,17 ]\,\,[ 16,15,11,12,14,10,13,1]$ 

$H4(8,8,8):=[0,21,22,19,20,23,18,9 ]\,\,[ 8,5,6,3,4,7,2,17 ]\,\,[ 16,13,14,11,12,15,10,1]$

$H5(8,8,8):=[0,21,22,19,20,23,10,17 ]\,\,[ 8,5,6,3,4,7,18,1 ]\,\,[ 16,13,14,11,12,15,2,9]$

$H6(8,8,8):=[0,21,22,19,20,23,9,18 ]\,\,[ 8,5,6,3,4,7,17,2 ]\,\,[ 16,13,14,11,12,15,1,10]$

$H7(8,8,8):=[0,20,19,22,21,18,9,23 ]\,\,[ 8,4,3,6,5,2,17,7 ]\,\,[ 16,12,11,14,13,10,1,15]$

$H0(6,6,6):=[1,20,19,6,21,18 ]\,\,[ 9,4,3,14,5,2 ]\,\,[ 17,12,11,22,13,10 ]\,\,[ 15,16,7,8,23,0]$

$H1(6,6,6):=[0,23,13,2,22,12 ]\,\,[ 8,7,21,10,6,20 ]\,\,[ 16,15,5,18,14,4 ]\,\,[ 11,17,3,9,19,1]$

$H2(6,6,6):=[1,20,19,0,23,18 ]\,\,[ 9,4,3,8,7,2 ]\,\,[ 17,12,11,16,15,10 ]\,\,[ 14,21,6,13,22,5]$

First Loupekine'snark= $[0,3,4,7,8,11 ]\,\,
[1,2,9,10,5,6]$

Second Loupekine's snark= $[0,9,10,7,4,1,2,11,8,5,6,3]$

The third $T_3$-snark $L_3:=[0, 4, 2, 1, 6, 8, 10, 5, 9, 11, 7, 3]$ 

($L_3$ is isomorphic to the $T_3$-snark $[0, 4, 8, 1, 5, 10] \,\, [2, 6, 9, 3, 7, 11]$)

Petersen graph = $[0,3,4,1,2,5]$

$Y$:=$[0,4,8,1,5,10 ]\,\,[2,6,12,3,7,14  ]\,\,[9,16,20,11,17,21 ]\,\,[ 13,18,22,15,19,23]$

%%%%%%%%%%%%%%%%%%%%%%%%%%%%%%%%%%%%%%%%%%%%%%%%%%%%%%%%%%%%%%%%%%%%%%%%%%%%%
%%%%%%%%%%%%%%%%%%%%%%%%%%%%%%%%%%%%%%%%%%%%%%%%%%%%%%%%%%%%%%%%%%%%%%%%%%%%%

\end{document}